\UseRawInputEncoding
\documentclass[12pt]{article}
\usepackage{amssymb}
\usepackage{amsthm}
\usepackage{amsthm,amsmath,indentfirst,amssymb}
\usepackage{pstricks}
\usepackage{algorithm} 
\usepackage{algorithmic} 
\usepackage{blindtext}
\usepackage{graphicx}
\usepackage{cite}
\usepackage{float}
\usepackage[colorlinks,
linkcolor=blue,
anchorcolor=blue,
citecolor=blue
]{hyperref}
\usepackage[scriptsize,nooneline]{subfigure}
\usepackage{epstopdf}
{}
{}
{}

\begin{document}

\title{The Boundedness of Fractional Integral Operators in Local and Global Mixed Morrey-type Spaces
\thanks{The research was supported by the National Natural Science Fundation of China(12061069).}
}
\author{Houkun Zhang,\quad Jiang Zhou\thanks{ Corresponding author E-mail address: zhoujiang@xju.edu.cn.}
\\
\small  College of Mathematics and System Sciences, Xinjiang University, Urumqi 830046\\
\small People's Republic of China
}
\date{}
\maketitle
{\bf Abstract:} In this paper, we introduced the local and global mixed Morrey-type spaces, and some properties of these spaces are also studied. After that,  the necessary conditions of the boundedness of fractional integral operators $I_{\alpha}$ are studied respectively in mixed-norm Lebesgue spaces and the local mixed Morrey-type spaces. Last but not least, the boundedness of $I_{\alpha}$ is obtained by Hardy operators' boundedness in weighted Lebesgue spaces. As corollaries, we acquire the boundedness of $I_{\alpha}$ in mixed Morrey spaces(\hyperlink{Corollary 5.6}{Corollary 5.6}) and mixed Lebesgue spaces(\hyperlink{Corollary 5.7}{Corollary 5.7}) respectively. Particularly, \hyperlink{Corollary 5.7}{Corollary 5.7} obtains a sufficient and necessary condition.
\par
{\bf Keywords:} Local and global mixed Morrey-type spaces; Fractional integral operators; Hardy operators

\maketitle

\section{Introduction}\label{sec1}

\par
If $E$ is a nonempty measurable subset on $\mathbb{R}^n$ and $f$ is a measurable function on $E$, then we put
$$\|f\|_{L_p(E)}:=\|f\chi_E\|_{L_p(\mathbb{R}^n)}=\left(\int_{\mathbb{R}^n}|f(y)|^p\chi_E dy\right)^{\frac{1}{p}}.$$

For $x\in \mathbb{R}^n$ and $r>0$, let $Q(x,r)$ denote the open cube centered at $x$ of length $2r$ and $Q(x,r)^{\complement}$ denote the set $\mathbb{R}^n\backslash Q(x,r).$

Let $\vec{p}=(p_1,p_2,\cdots,p_n),~\vec{q}=(q_1,q_2,\cdots,q_n),~\vec{s}=(s_1,s_2,\cdots,s_n)~\cdots$ are n-tuples and $0<p_i,q_i,s_i<\infty,~i=1,2,\cdots,n$. We define that  $\frac{1}{\vec{p}}=(\frac{1}{p_1},\frac{1}{p_2},\cdots,\frac{1}{p_n})$ and $\vec{p}<\vec{q}$ means $p_i<q_i$ holds for each $i$. In particular, $\vec{p}<r$ means $p_i<r$ hold for each $i$.

Let $f\in L^{loc}_1(\mathbb{R}^n)$. The fractional integral operators $I_\alpha$ were defined by
$$I_\alpha f(x)=\int_{\mathbb{R}^n}\frac{f(y)}{|x-y|^{n-\alpha}} \,dy,~~0<\alpha<n.$$
These operators $I_\alpha$ play an essential role in real and harmonic analysis \cite{1,2}.

In partial differential equations, Morrey spaces $\mathcal{M}_{p,\lambda}$ are widely used to investigate the local behavior of
solutions to elliptic and parabolic differential equations. In 1938, Morrey introduced these spaces \cite{3} and they were defined as following: For $0\le\lambda\le n,~1\le p\le\infty$, we say that $f\in\mathcal{M}_{p,\lambda}$ if $f\in L^{loc}_{p}(\mathbb{R}^n)$ and
$$\|f\|_{\mathcal{M}_{p,\lambda}}=\|f\|_{\mathcal{M}_{p,\lambda}(\mathbb{R}^n)}=\sup_{x\in\mathbb{R}^n,~r>0}r^{-\frac{\lambda}{p}}\|f\|_{L_p(Q(x,r))}<\infty.$$
It is obvious that if $\lambda=0$ then $\mathcal{M}_{p,0}=L_p(\mathbb{R}^n)$; if $\lambda=n$, then $\mathcal{M}_{p,n}=L_{\infty}(\mathbb{R}^n)$; if $\lambda<0$ or $\lambda>n$, then $\mathcal{M}_{p,\lambda}=\Theta$, where $\Theta$ is the set of all functions almost everywhere equivalent to 0 on $\mathbb{R}^n$.

For nearly two decades, due to the more precise structure of mixed-norm function spaces than the corresponding classical function spaces, the mixed-norm function spaces are widely used in the partial differential equations \cite{20,21,22,23}.

In 1961, the mixed Lebesgue spaces $L_{\vec{p}}(\mathbb{R}^n)$ were studied by Benedek and Panzone \cite{15}. These spaces were
natural generalizations of the classical Lebesgue spaces $L_p$. The definition was stated as following: Let $f$ is a measurable function on $\mathbb{R}^n$ and  $0<\vec{p}\le\infty$; we say that $f$ belongs to the mixed Lebesgue spaces $L_{\vec{p}}(\mathbb{R}^n)$, if the number obtained after taking successively the $p_1$-norm for $x_1$, the $p_2$-norm for $x_2$, $\cdots$ , the $p_n$-norm for $x_n$, and in that order, is finite;
the number obtained will be denoted by $\left\|f\right\|_{L_{\vec{p}}(\mathbb{R}^{n})}$ or $\|\cdots\|~\|f\|_{L_{p_1}(\mathbb{R})}~\|_{L_{p_2}(\mathbb{R})}\cdots\|_{L_{p_n}(\mathbb{R})}$.

When $\vec{p}<\infty$, we can write $\left\|f\right\|_{L_{\vec{p}}(\mathbb{R}^{n})}$ as following,
$$\left\|f\right\|_{L_{\vec{p}}(\mathbb{R}^n)}=\left(\int_{\mathbb{R}}\cdots\left(\int_{\mathbb{R}}\left|f(x)\right|^{p_1}\,dx_1\right)
^{\frac{p_2}{p_1}}\cdots\,dx_n\right)^{\frac{1}{p_n}}<\infty.$$
If $p_1=p_2=\cdots=p_n=p$, then $L_{\vec{p}}(\mathbb{R}^n)$ are reduced to classical Lebesgue spaces $L_p$ and
$$\left\|f\right\|_{L_{\vec{p}}(\mathbb{R}^n)}=\left(\int_{\mathbb{R}^n}\left|f(x)\right|^{p} dx\right)^{\frac{1}{p}}.$$
Particularly, if $p_1=p_2=\cdots=p_n=\infty$, then the result is similar: $L_{\vec{p}}=L_{\infty}$.

After that, a host of function spaces with mixed norm were introduced, such as mixed-norm Lorentz spaces \cite{24}, mixed-norm Lorentz-Marcinkiewicz spaces \cite{27}, mixed-norm Orlicz spaces \cite{25}, anisotropic mixed-norm Hardy spaces \cite{26}, mixed-norm Triebel-Lizorkin spaces \cite{28} and weak mixed-norm Lebesgue spaces \cite{29}. More information can be found in \cite{33}.

Moreover, combining mixed Lebesgue spaces and Morrey spaces, Nogayama, in 2019, introduced mixed Morrey spaces \cite{4,5} and these spaces were defined as following: let $\vec{q}=(q_1,q_2,\cdots,q_n)\in (0,\infty]^n$ and $p\in(0,\infty]$ satisfy
$$\sum_{j=1}^{n}\frac{1}{q_j}\ge\frac{n}{p};$$
the mixed Morrey spaces $\mathcal{M}_{\vec{q}}^{p}(\mathbb{R}^n)$ were defined to be the set of all measurable functions $f$ such that their quasi-norms
$$\|f\|_{\mathcal{M}_{\vec{q}}^p(\mathbb{R}^n)}:=\sup\bigg\{|Q|^{\frac{1}{p}-\frac{1}{n}(\sum_{j=1}^n\frac{1}{q_j})}\|f\chi_Q\|_{L_{\vec{q}}(\mathbb{R}^n)}:Q~is~ a~cube~in~\mathbb{R}^n\bigg\}$$
are finite. It is obvious that if $q_1=q_2=\cdots=q_n=q$, then $\mathcal{M}_{\vec{q}}^p=\mathcal{M}_q^p=\mathcal{M}_{q,\lambda}$ and if $\frac{1}{p}=\frac{1}{n}\sum_{j=1}^n\frac{1}{q_j}$, then $\mathcal{M}_{\vec{q}}^p=L_{\vec{q}}.$

In particular, before these definitions of mixed Morrey spaces, another definitions were given by \cite{30} in 2017. Their norms were defined by taking successively the $\mathcal{M}_{q_1}^{p_1}$-norm for $x_1$, the $\mathcal{M}_{q_2}^{p_2}$-norm for $x_2$. It is easy to prove that
$$\|f\|_{\mathcal{M}_{\vec{q}}^p(\mathbb{R}^2)}\le \left\|\left\|f\right\|_{\mathcal{M}_{q_1}^{p_1}(\mathbb{R})}\right\|_{\mathcal{M}_{q_2}^{p_2}(\mathbb{R})},$$
if $\frac{1}{p}=\frac{1}{2}\sum_{i=1}^{2}\frac{1}{p_i}$.

In 1981, Adams introduced a variant of Morrey-typle spaces $\mathcal{M}_{p\theta,\lambda}$ \cite{19}.
In 2004, Burenkov and Guliyev developed these spaces and defined local Morrey-type spaces $LM_{p\theta,\omega}$ and global Morrey-type spaces $GM_{p\theta,\omega}$ with their some properties were studied \cite{6}. They are defined as follows: let $0<p,\theta\le\infty$ and let $\omega$ be a non-negative measurable function on $(0,\infty)$. We say $f$ belong to the local Morrey-type spaces $LM_{p\theta,\omega}$, if the quasi-norms
$$\|f\|_{LM_{p\theta,\omega}}=\|f\|_{LM_{p\theta,\omega}(\mathbb{R}^n)}=\|\omega(r)\|f\|_{L_p(Q(0,r))}\|_{L_{\theta}(0,\infty)}.$$
are finite and $f$ belong to the global Morrey-type spaces $GM_{p\theta,\omega}$, if the quasi-norms
$$\|f\|_{GM_{p\theta,\omega}}=\sup_{x\in\mathbb{R}^n}\|f(x+\cdot)\|_{LM_{p\theta,\omega}}.$$
are finite.

Note that $GM_{p\theta,r^{-\lambda/p}}=\mathcal{M}_{p\theta,\lambda},~0\le\lambda\le n$ and $GM_{p\infty,r^{-\lambda/p}}=\mathcal{M}_{p,\lambda},~0\le\lambda\le n$.

This paper will give the definitions of local and global mixed Morrey-type spaces and their some properties.

In 1955, Plessis studied the boundedness of fractional integral operators \cite{7} in classical Lebesgue spaces. Adams in 1975 investigated boundedness of the fractional integral operators in Morrey spaces \cite{8}. If $\omega$, which is a positive measurable function defined on $(0,\infty)$, replace the power function $r^{-\lambda/p}$ in the definitions of $\mathcal{M}_{p,\lambda}$, then these become the Morrey-type spaces $\mathcal{M}_{p,\omega}$. Mizuhara \cite{9}, Nakai \cite{10}, and Guliyev \cite{11} generalized Adams' result respectively. The boundedness of fractional integral operators in local Morrey-type spaces was also studied in \cite{12,13,14}.

In \cite{15}, Benedek and Panzone proved the boundedness of fractional integral operators in mixed lebesgue spaces and it was stated as follows: Let $X=E^n$, and $\vec{\alpha}=(\alpha_1,\alpha_2,\cdots,\alpha_n)$ an n-tuple of real numbers, $0<\alpha_i<1$. If $\vec{p}$ and $\vec{q}$ are such that $\frac{1}{\vec{p}}-\frac{1}{\vec{q}}=\vec{\alpha},~1<\vec{p}<\frac{1}{\vec{\alpha}}$, then
$$\|I_\beta f\|_{L_{\vec{q}}}\le C\|f\|_{L_{\vec{p}}}$$
holds for every $f\in L_{\vec{p}}$, where $\beta=\sum_{i=1}^{n}\alpha_i$ and $C= C(\vec{\alpha},\vec{p})$.

In 1974, Adams and David studied the boundedness of the fractional integral operators $I_{\alpha}$ in the mixed lebesgue spaces $L_{\vec{p}}(\mathbb{R}^2)$ \cite{31}. Their results are stated as the following: Let $\vec{p}=(p_1,p_2),~\vec{q}=(q_1,q_2),~1\le p_1\le q_1\le\infty,~1<p_2<q_2<\infty$, and $\alpha=\sum_{i=1}^2\frac{1}{p_i}-\sum_{i=1}^2\frac{1}{q_i}$, we have
$$\|I_\alpha f\|_{L_{\vec{q}}}\le C\|f\|_{L_{\vec{p}}};$$
let $\vec{p}=(p_1,p_2),~\vec{q}=(q_1,q_2),~1<p_1<q_1<\infty,~1<p_2=q_2<\infty$, and $\alpha=\frac{1}{p_1}-\frac{1}{q_1}$, then
$$\|I_\alpha f\|_{L_{\vec{q}}}\le C\|f\|_{L_{\vec{p}}}.$$
It is unlucky that we fail to prove that the second result is also right when $n>2$ and $p_2,q_2$ are replaced by $\mathbf{p'}=(p_{l+1},\cdots,p_n),\mathbf{q'}=(q_{l+1},\cdots,q_n)(n-l\ge 2)$ respectively because of the application of Fubini's theorem.

In \cite{4}, Toru Nogayama proved the boundedness of fractional integral operators in mixed Morrey spaces and it was stated as following: Let $0<\alpha<n,~1<\vec{q}, \vec{s}<\infty$ and $0<p,r<\infty$; assume that $\frac{n}{p}\le\sum_{j=1}^n\frac{1}{q_j}$ and $\frac{n}{r}\le\sum_{j=1}^n\frac{1}{s_j}$; if
$$\frac{1}{r}=\frac{1}{p}-\frac{\alpha}{n},~\frac{\vec{q}}{p}=\frac{\vec{s}}{r};$$
then, for $f\in\mathcal{M}_{\vec{q}}^p(\mathbb{R}^n)$,
$$\|I_{\alpha}f\|_{\mathcal{M}_{\vec{s}}^r(\mathbb{R}^n)}\le C\|f\|_{\mathcal{M}_{\vec{q}}^p(\mathbb{R}^n)}.$$
It is obvious that this result can not be reduced to the result of the above, even if $\frac{1}{p}=\frac{1}{n}\sum_{j=1}^n\frac{1}{q_j},~\frac{1}{r}=\frac{1}{n}\sum_{j=1}^n\frac{1}{s_j}$.

Let $\vec{p}_1=(p_{11},p_{12},\cdots,p_{1n}),~\vec{p}_2=(p_{21},p_{22},\cdots,p_{2n})$, and $0<p_{1i},p_{2i}<\infty, i-1,2,\cdots,n$.
Moreover, we define that $f\in L^{loc}_{\vec{p}}$ means that $f\chi_E\in L_{\vec{p}}$, if $E$ is compact set. If $E$ is a nonempty measurable subset on $\mathbb{R}^n$ and $f$ is a measurable function on $E$, then we put
$$\|f\|_{L_{\vec{p}}(E)}:=\|f\chi_E\|_{L_{\vec{p}}(\mathbb{R}^n)}.$$
By $A\lesssim B(B\gtrsim A)$, we denote that $A\le CB$ where $C>0$ depends on unessential parameters, and $A\thicksim B$ means that $A\lesssim B$ and $A\gtrsim B$. We define that $t_+=t$ if $t>0$ and $t_+=0$ if $t<0$.

The paper is organized as follows. We start with defining the local and global mixed Morrey-type spaces, and some properties of these spaces are also studied in Section 2. In Section 3,  We will give the necessary conditions of the boundedness of fractional integral operators $I_{\alpha}$ in mixed-norm Lebesgue spaces and the local mixed Morrey-type spaces. In order to acquire relationship between fractional integral operators $I_{\alpha}$ and Hardy operators $H$, $L_{\vec{p}}\,$-estimates of $I_{\alpha}$ over the cube $Q(x,r)$ is investigated in the Section 4. In the last section, according to results of Section 4, the inequality $\|I_{\alpha}f\|_{LM_{\vec{p}_2\theta_2,\omega_2}} \lesssim\|Hg_{\vec{p}_1}\|_{L{\theta_2,\nu_2}(0,,\infty)}$ is proved. As corollaries, we acquire the boundedness of fractional integral operators in global mixed Morrey-type spaces, the boundedness of fractional maximal operators in local and global mixed Morrey-type spaces, and the other two new results of the boundedness of fractional integral operators in mixed Lebesgue spaces and mixed Morrey spaces.

\section{Local and Global Mixed Morrey-type Spaces}\label{sec1}

\par
In this section, we will give the definitions of local and global mixed Morrey-type spaces.

From \cite{15}, we know H\"{o}lder's inequality and Minkowski's inequality in mixed Lebesgue spaces. In \cite{4}, Fatou's property for $L_{\vec{p}}$ was also given.

\textbf{Lemma 2.1} (H\"{o}lder's inequality for $L_{\vec{p}}$). Let $1\le\vec{p}\le\infty$ and $\frac{1}{\vec{p}}+\frac{1}{\vec{p}\,'}=1$. Then for any $f\in L_{\vec{p}}$ and $g\in L_{\vec{p}\,'}$,
$$\int_{\mathbb{R}^n}f(x)g(x)\,dx\le\|f\|_{L_{\vec{p}}}\|g\|_{L_{\vec{p}\,'}}$$
holds.

\textbf{Remark 2.2.} In particular, if $0<\vec{r},\vec{p},\vec{q}<\infty$, $\frac{1}{\vec{p}}+\frac{1}{\vec{q}}=\frac{1}{\vec{r}}$, $f\in L_{\vec{p}}$ and $g\in L_{\vec{q}}$, then
$$\|fg\|_{L_{\vec{r}}}\le\|f\|_{L_{\vec{p}}}\|g\|_{L_{\vec{q}}}$$
also holds.

They are easily proved by successive applications of classical H\"{o}lder's inequality.

\textbf{Lemma 2.3} (Minkowski's inequality for $L_{\vec{p}}$). Let $1\le\vec{p}\le\infty$. If $f,g\in L_{\vec{p}}$, then
$$\|f+g\|_{L_{\vec{p}}}\le\|f\|_{L_{\vec{p}}}+\|g\|_{L_{\vec{p}}}.$$

It is easy to prove it by successive applications of classical Minkowski's inequality.

\textbf{Lemma 2.4} (Fatou's property for $L_{\vec{p}}$). Let $0<\vec{p}\le\infty$. Let $\{f_i\}_{i=1}^{\infty}$ be a sequence of non-negative measurable functions on $\mathbb{R}^n$. Then
$$\left\|\liminf_{i\rightarrow\infty}f_i\right\|\le\liminf_{i\rightarrow\infty}\left\|f_i\right\|.$$

\textbf{Remark 2.5.} By classical Fatou's property, it also can be proved that
$$\limsup_{i\rightarrow\infty}\left\|f_i\right\|\le\left\|\limsup_{i\rightarrow\infty}f_i\right\|.$$

We state the Lebesgue differential theorem in the setting of mixed-norm Lebesgue spaces as the following lemma.

\hypertarget{Lemma 2.6}{\textbf{Lemma 2.6.}} Let $f\in L^{loc}_{\vec{p}}$ and $0<\vec{p}<\infty$, then
$$\lim_{r\rightarrow 0}\|\chi_{Q(x,r)}\|_{L_{\vec{p}}(\mathbb{R}^n)}^{-1}\|f\|_{L_{\vec{p}}(Q(x,r))}=|f(x)|~~a.e.~x\in\mathbb{R}^n.$$

\textbf{Proof:} Let $\vec{p}=(p_1,p_2,\cdots,p_n)$. Without less of generality, we assume $n=2$.

It is easy to know that
$$\|\chi_{Q(x,r)}\|_{L_{\vec{p}}}=|2r|^{\sum_{i=1}^{2}\frac{1}{p_i}},$$
where $x=(x_1,x_2).$
Therefore, by the classical Lebesgue differential theorem, we get that for any $\epsilon>0$, there exist $\delta>0$, such that $r<\delta$ and
\begin{eqnarray*}
\|\chi_{Q(x,r)}\|_{L_{\vec{p}}(\mathbb{R}^n)}^{-1}\|f\|_{L_{\vec{p}}(Q(x,r))}
&=&|2r|^{-\sum_{i=1}^{2}\frac{1}{p_i}}\bigg(\int_{\mathbb{R}}\Big(\int_{\mathbb{R}}|f(y_1,y_2)|^{p_1}\chi_{Q(x,r)}(y_1,y_2)\,dy_1\Big)^{\frac{p_2}{p_1}}\,dy_2\bigg)^{\frac{1}{p_2}}\\
&=&\bigg(\frac{1}{2r}\int_{I_2}\Big(\frac{1}{2r}\int_{I_1}|f(y_1,y_2)|^{p_1}\,dy_1\Big)^{\frac{p_2}{p_1}}\,dy_2\bigg)^{\frac{1}{p_2}}\\
&\le&\left(\left(|f(x_1,x_2)|^{p_1}+\epsilon\right)^{\frac{p_2}{p_1}}+\epsilon\right)^{\frac{1}{p_2}}~~a.e.~x\in\mathbb{R}^n\\
&\lesssim&|f(x_1,x_2)|+\epsilon^{\frac{1}{p_1}}+\epsilon^{\frac{1}{p_2}}
\end{eqnarray*}
where $I_1=(x_1-r,x_1+r),~I_2=(x_2-r,x_2+r)$. By the definition of limit,
$$\lim_{r\rightarrow 0}\|\chi_{Q(x,r)}\|_{L_{\vec{p}}(\mathbb{R}^n)}^{-1}\|f\|_{L_{\vec{p}}(Q(x,r))}=|f(x)|~~a.e.~x\in\mathbb{R}^n.$$

The proof is complete.    $~~~~\blacksquare$

Based on the definition of local and global Morrey-type spaces and mixed Morrey spaces, the definitions of local and global mixed Morrey-type spaces are introduced.

\textbf{Definition 2.7.} Let $0<\vec{p},\theta\le\infty$ and let $\omega$ be a non-negative measurable function on $(0,\infty)$. We denote the local mixed Morrey-type spaces and the global mixed Morrey-type spaces by $LM_{\vec{p}\theta,\omega}$, $GM_{\vec{p}\theta,\omega}$ respectively. For any functions $f\in L^{loc}_{\vec{p}}(\mathbb{R}^n)$, we say $f\in LM_{\vec{p}\theta,\omega}$ when the quasi-norms
$$\|f\|_{LM_{\vec{p}\theta,\omega}}=\|f\|_{LM_{\vec{p}\theta,\omega}(\mathbb{R}^n)}=\|\omega(r)\|f\|_{L_{\vec{p}}(Q(0,r))}\|_{L_{\theta}(0,\infty)} \le\infty;$$
we say $f\in GM_{\vec{p}\theta,\omega}$ when the quasi-norms
$$\|f\|_{GM_{\vec{p}\theta,\omega}}=\sup_{x\in\mathbb{R}^n}\|f(x+\cdot)\|_{LM_{\vec{p}\theta,\omega}}\le\infty.$$
In particular, we say $f\in LM_{\vec{p}\theta,\omega}^{[x]}$, if
$$\|f(x+\cdot)\|_{LM_{\vec{p}\theta,\omega}}=\|\omega(r)\|f\|_{L_{\vec{p}}(Q(x,r))}\|_{L_{\theta}(0,\infty)}<\infty.$$

\textbf{Remark 2.8.} Note that if $\omega(r)=r^{\frac{n}{q}-\sum_{j=1}^{n}\frac{1}{p_j}},~\theta=\infty$, then $GM_{\vec{p}\theta,\omega}=\mathcal{M}_{\vec{p}}^{q}$.

Next, some properties of local and global mixed Morrey-type spaces are given.

\textbf{Theorem 2.9.} Let $0<\vec{p},\theta\le\infty$ and let $\omega$ be a non-negative measurable function on $(0,\infty)$.\\
(1) If for all $t>0$
$$\|\omega(r)\|_{L_{\theta}(t,\infty)}=\infty,\hypertarget{2.1}{\eqno{(2.1)}}$$
then $LM_{\vec{p}\theta,\omega}=GM_{\vec{p}\theta,\omega}=\Theta.$\\
(2) If for all $t>0$
$$\|\omega(r)r^{\sum_{j=1}^{n}\frac{1}{p_i}}\|_{L_{\theta}(0,t)}=\infty,\hypertarget{2.2}{\eqno{(2.2)}}$$
then $f(0)=0$ for all $f\in LM_{\vec{p}\theta,\omega}$ continuous at 0 and $GM_{\vec{p}\theta,\omega}=\Theta$ for all $0<\vec{p}<\infty$.

\textbf{Proof:} (1) Suppose that (\hyperlink{2.1}{2.1}) holds true for all $t\in(0,\infty)$ and $f$ is not the element zero. Then there is $t_0$ such that $A=\|f\|_{L_{\vec{p}}(Q(0,t_0))}>0$. Hence,
$$\|f\|_{GM_{\vec{p}\theta,\omega}}\ge\|f\|_{LM_{\vec{p}\theta,\omega}}\ge\|\omega(r)\|f\|_{L_{\vec{p}}(Q(0,r))}\|_{L_{\theta}(t_0,\infty)}\ge A\|\omega(r)\|_{L_{\theta}(t_0,\infty)}=\infty.$$
Therefore, $\|f\|_{GM_{\vec{p}\theta,\omega}}=\|f\|_{LM_{\vec{p}\theta,\omega}}=\infty.$\\
(2) Suppose that (\hyperlink{2.2}{2.2}) holds true for all $t\in(0,\infty)$ and $f\in LM_{\vec{p}\theta,\omega}$. By \hyperlink{Lemma 2.6}{Lemma 2.6} there  exists
$$\lim_{r\rightarrow 0}r^{-\sum_{j=1}^{n}\frac{1}{p_i}}\|f\|_{L_{\vec{p}}(Q(0,r))}=B=C|f(0)|~~a.e.~x\in\mathbb{R}^n.$$
where $C=2^{\sum_{j=1}^{n}\frac{1}{p_i}}$. If $B>0$, then there exists $t_0>0$ suc that
$$r^{-\sum_{j=1}^{n}\frac{1}{p_i}}\|f\|_{L_{\vec{p}}(Q(0,r))}\ge\frac{B}{2}~~a.e.~x\in\mathbb{R}^n$$
holds for all $0<r\le t_0$. Consequently,
$$\|f\|_{LM_{\vec{p}\theta,\omega}}\ge\|\omega(r)\|f\|_{L_{\vec{p}}(Q(0,r))}\|_{L_{\theta}(0,t_0)}\ge\frac{B}{2}\|\omega(r)r^{\sum_{j=1}^{n}\frac{1}{p_i}}\|_{L_{\theta}(0,t)}=\infty~~a.e.~x\in\mathbb{R}^n.$$
Here, $\|f\|_{LM_{\vec{p}\theta,\omega}}=\infty$, so $f\notin LM_{\vec{p}\theta,\omega}$ and we have arrived at a contradiction.

Next let $0<\vec{p}<\infty,~f\in GM_{\vec{p}\theta,\omega}$. By lemma2.1, we know that
$$\lim_{r\rightarrow 0}r^{-\sum_{j=1}^{n}\frac{1}{p_i}}\|f\|_{L_{\vec{p}}(Q(x,r))}=C|f(x)|~~a.e.~x\in\mathbb{R}^n,$$
where $C=2^{\sum_{j=1}^{n}\frac{1}{p_i}}$. By the above argument, $f(x)=0$ a.e. $x\in\mathbb{R}^n$. $~~~~\blacksquare$

Due to the above argument, the sets $\Omega_{\theta}$ and $\Omega_{\vec{p},\theta}$ are defined as following.

\textbf{Definition 2.10.} Let $0<\vec{p},\theta\le\infty$. We denote by $\Omega_{\theta}$ the set of all functions $\omega$ which are non-negative, measurable on $(0,\infty)$, not equivalent to 0, and such that for some $t>0$
$$\|\omega\|_{L_{\theta}(t,\infty)}<\infty.$$
Moreover, we denote by $\Omega_{\vec{p},\theta}$ the set of all functions $\omega$ which are non-negative, measurable on $(0,\infty)$, not equivalent to 0, and such that for some $t_1,t_2>0$
$$\|\omega\|_{L_{\theta}(t_1,\infty)}<\infty,~~~~\|\omega(r)r^{\sum_{j=1}^{n}\frac{1}{p_i}}\|_{L_{\theta}(0,t_2)}=\infty.$$

The following Theorem tells us that $LM_{\vec{p}\theta,\omega}$ and $GM_{\vec{p}\theta,\omega}$ are complete, when $1\le\vec{p},\theta<\infty$.

\textbf{Theorem 2.11.} (1) Let $1\le\vec{p},\theta<\infty$, and $\omega\in\Omega_{\theta}$. Suppose $f_n\in LM_{\vec{p}\theta,\omega}$, for all $n\in \mathbb{N}_+$ and
$$\sum_{n=1}^{\infty}\|f_n\|_{LM_{\vec{p}\theta,\omega}}<\infty.$$
Then $\sum_{n=1}^{\infty}f_n$ exist. If $f=\sum_{n=1}^{\infty}f_n$, then $f\in LM_{\vec{p}\theta,\omega}$ and
$$\|f\|_{LM_{\vec{p}\theta,\omega}}\le\sum_{n=1}^{\infty}\|f_n\|_{LM_{\vec{p}\theta,\omega}}.$$
Hence, $LM_{\vec{p}\theta,\omega}$ is compte.\\
(2) Let $1\le\vec{p},\theta<\infty$, and $\omega\in\Omega_{\vec{p}\theta}$. Suppose $f_n\in GM_{\vec{p}\theta,\omega}$, for all $n\in \mathbb{N}_+$ and
$$\sum_{n=1}^{\infty}\|f_n\|_{GM_{\vec{p}\theta,\omega}}<\infty.$$
Then $\sum_{n=1}^{\infty}f_n$ exist. If $f=\sum_{n=1}^{\infty}f_n$, then $f\in GM_{\vec{p}\theta,\omega}$ and
$$\|f\|_{GM_{\vec{p}\theta,\omega}}\le\sum_{n=1}^{\infty}\|f_n\|_{GM_{\vec{p}\theta,\omega}}.$$
Hence, $GM_{\vec{p}\theta,\omega}$ is compte.

\textbf{Proof:} (1) It is easy to see that for any $R>0$
$$\|\omega\|_{L_{\theta}(R,\infty)}\|f\|_{L_{\vec{p}}(Q(0,R))}\le\|f\|_{LM_{\vec{p}\theta,\omega}}.$$
Thus
$$\sum_{n=1}^{\infty}\|f_n\|_{L_{\vec{p}}(Q(0,R))}\le C\sum_{n=1}^{\infty}\|f_n\|_{LM_{\vec{p}\theta,\omega}}.$$
Since $L_{\vec{p}}(Q(0,r))$ is complete \cite{15}, then $\sum_{n=1}^{\infty}f_n$ converges a.e. to some $f\in L^{loc}_{\vec{p}}(\mathbb{R}^n)$
$$\sum_{n=1}^{\infty}f_n=f$$
and
$$\|f\|_{L_{\vec{p}}(Q(0,R))}\le\sum_{n=1}^{\infty}\|f_n\|_{L_{\vec{p}}(Q(0,R))}.$$
Therefor,
\begin{eqnarray*}
\|f\|_{LM_{\vec{p}\theta,\omega}}&=&\|\omega(r)\|f\|_{L_{\vec{p}}(Q(0,r))}\|_{L_{\theta}(0,\infty)}\\
&\le&\|\omega(r)\sum_{n=1}^{\infty}\|f_n\|_{L_{\vec{p}}(Q(0,r))}\|_{L_{\theta}(0,\infty)}\\
&\le&\sum_{n=1}^{\infty}\|\omega(r)\|f_n\|_{L_{\vec{p}}(Q(0,r))}\|_{L_{\theta}(0,\infty)}\\
&=&  \sum_{n=1}^{\infty}\|f_n\|_{LM_{\vec{p}\theta,\omega}}
\end{eqnarray*}

Now, the complete of $LM_{\vec{p}\theta,\omega}$ will be proved. let $\{f_n\}$ is Cauchy sequence in $LM_{\vec{p}\theta,\omega}$. Without less of generality, assume $\{f_n\}$ satisfies that
$$\sum_{n=1}^{\infty}\|f_n-f_{n-1}\|_{LM_{\vec{p}\theta,\omega}}<\infty,$$
where $f_0=0$.

Thanks to the above argument, there exist $f=\sum_{n=1}^{\infty}(f_n-f_{n-1})=\lim_{n\rightarrow\infty}f_n$ and $f\in LM_{\vec{p}\theta,\omega}$. It is easy via Fatou's property to obtain that
$$\limsup_{n\rightarrow\infty}\|f-f_n\|_{LM_{\vec{p}\theta,\omega}}\le\|\limsup_{n\rightarrow\infty}|f-f_n|\|_{LM_{\vec{p}\theta,\omega}}=0$$
Hence, the complete of $LM_{\vec{p}\theta,\omega}$ is proved.

(2) We can prove (2) through the same method as (1), so we omit the proof. $~~~~\blacksquare$

\textbf{Definition 2.12.}\cite{34} The intersection of a family function space $\{X_{\alpha}\}_{\alpha\in A}$ is a Banach space $X$ such that

(a) $X\hookrightarrow X_{\alpha},~\alpha\in A$;

(b) if for a certain Banach Space $Y$ we have
$$Y\hookrightarrow X_{\alpha},~\alpha\in A$$
then $Y\hookrightarrow X$.

It is easy to prove the following prosperity.

\textbf{Prosperity 2.13.} $GM_{\vec{p}\theta,\omega}$ is the interaction of $\{LM_{\vec{p}\theta,\omega}^{[x]}\}_{x\in \mathbb{R}^n}$.

\section{the necessary condition of the boundedness of $I_{\alpha}$}\label{sec1}

\par
In \cite{15,31}, the boundedness of fractional integral operators is proved in mixed Lebesgue spaces. Next, a sufficient and necessary condition of fractional integral operators' boundedness will be given in mixed norm Lebesgue spaces.

\hypertarget{Lemma 3.1}{\textbf{Lemma 3.1.}} Let $0<\alpha<n,~1<\vec{p}<\vec{q}<\infty$. Then
$$\|I_{\alpha}f\|_{L_{\vec{q}}}\lesssim\|f\|_{L_{\vec{p}}}\eqno{(3.1)}$$
if and only if
$$\alpha=\sum_{j=1}^{n}\frac{1}{p_j}-\sum_{j=1}^{n}\frac{1}{q_j}.$$

\textbf{Proof:} Without less of generality, let $n=2$. Duce to $|y-z|=|(y_1,y_2)-(z_1,z_2)|\ge|y_i-z_i|,~i=1,2$, it is easy to acquire that
$$|y-z|^{\alpha-2}\le\prod_{i=1}^{2}|y_i-z_i|^{\alpha_i-1},$$
where $\alpha_i=\frac{1}{p_i}-\frac{1}{q_i},~\alpha_i\in(0,1)$. Duce to Minkowski's inequality and boundedness of $I_{\alpha}$ in classical Lebesgue spaces, It is easy to acquire that
\begin{eqnarray*}
\|I_{\alpha}f\|_{L_{\vec{q}(\mathbb{R}^2)}}&=&\bigg(\int_{\mathbb{R}}\bigg(\int_{\mathbb{R}}\bigg|\int_{\mathbb{R}}\int_{\mathbb{R}}\frac{f(z_1,z_2)}{|x-z|^{2-\alpha}}\,dz_1\,dz_2\bigg|^{q_1}\,dy_1\bigg)^{\frac{q_2}{q_1}}\,dy_2\bigg)^{\frac{1}{q_2}}\\
&=&\bigg(\int_{\mathbb{R}}\bigg(\int_{\mathbb{R}}\frac{1}{|x_2-z_2|^{1-{\alpha}_2}}\bigg(\int_{\mathbb{R}}\bigg|\int_{\mathbb{R}}\frac{f(z_1,z_2)}{|x_1-z_1|^{1-{\alpha}_1}}\,dz_1\bigg|^{q_1}\,dy_1\bigg)^{\frac{1}{q_1}}\,dz_2\bigg)^{q_2}\,dy_2\bigg)^{\frac{1}{q_2}}\\
&=&\bigg(\int_{\mathbb{R}}\bigg(\int_{\mathbb{R}}\frac{\|f(\cdot,z_2)\|_{L_{p_1}(\mathbb{R})}}{|x_2-z_2|^{1-{\alpha}_2}}\,dz_2\bigg)^{q_2}\,dy_2\bigg)^{\frac{1}{q_2}}\\
&=&\|f\|_{L_{\vec{p}}(\mathbb{R}^2)}.
\end{eqnarray*}
Hence, sufficiency is proved.

Let $\delta_tf(x)=f(tx)$. Then
$$\delta_{t^{-1}}I_{\alpha}\delta_t=t^{-\alpha}I_{\alpha},~\|\delta_tf\|_{L_{\vec{p}}(\mathbb{R}^2)}=t^{-\sum_{i=1}^2\frac{1}{p_i}}\|f\|_{L_{\vec{p}}(\mathbb{R}^2)}.$$
$$\|\delta_{t^{-1}}I_{\alpha}f\|_{L_{\vec{q}}(\mathbb{R}^2)}=t^{\sum_{i=1}^2\frac{1}{p_i}}\|I_{\alpha}f\|_{L_{\vec{q}}(\mathbb{R}^2)}.$$
Assume (3.1) is satisfied. Then
\begin{eqnarray*}
\|I_{\alpha}f\|_{L_{\vec{q}}(\mathbb{R}^2)}&=&t^{\alpha}\|\delta_{t^{-1}}I_{\alpha}\delta_tf\|_{L_{\vec{q}}(\mathbb{R}^2)}=t^{\alpha+\sum_{i=1}^2\frac{1}{q_i}}\|I_{\alpha}\delta_tf\|_{L_{\vec{q}}(\mathbb{R}^2)}\\
&\le&Ct^{\alpha+\sum_{i=1}^2\frac{1}{q_i}}\|\delta_tf\|_{L_{\vec{p}}(\mathbb{R}^2)}=Ct^{\alpha+\sum_{i=1}^2\frac{1}{q_i}-\sum_{i=1}^2\frac{1}{p_i}}\|f\|_{L_{\vec{p}}(\mathbb{R}^2)}.
\end{eqnarray*}
Therefor, $\alpha=\sum_{i=1}^2\frac{1}{p_i}-\sum_{i=1}^2\frac{1}{q_i}$ is obtained. $~~~~\blacksquare$

\hypertarget{Remark 3.2}{\textbf{Remark 3.2.}} By translation invariance of $I_{\alpha}$, we can prove that $\vec{p}\le\vec{q}$ is a necessary condition for the boundedness of fractional integral operators in mixed Lebesgue spaces.

Let $\tau^i_{h}f(x)=f(x_1,\cdots,x_{i-1},x_i+h,x_{i+1},\cdots,x_n)(i=1,2,\cdots,n)$ and $\|I_{\alpha}f\|_{L_{\vec{q}}}\le C\|f\|_{L_{\vec{p}}}$. By translation invariant of $I_{\alpha}$,
$$\|I_{\alpha}f+\tau^i_{h}I_{\alpha}f\|_{L_{\vec{q}}}=\|I_{\alpha}(f+\tau^i_{h}f)\|_{L_{\vec{q}}}\le C\|f+\tau^i_{h}f\|_{L_{\vec{p}}}.$$
Duce to Lemma 2.3 of \cite{32}, let $h\rightarrow\infty$, then
$$\|I_{\alpha}f\|_{L_{\vec{q}}}\le 2^{\frac{1}{p_i}-\frac{1}{q_i}}C\|f\|_{L_{\vec{p}}}.$$
Therefore, $\vec{p}\le\vec{q}$.

The following theorem states the necessity of $I_{\alpha}$ in local mixed Morrey-type spaces.

\textbf{Theorem 3.3.} (1) Let $1\le\vec{p_1}\le\infty,~0<\vec{p_2}\le\infty,~0<\alpha<n,~0<\theta_1,\theta_2<\infty,~\omega_1\in\Omega_{\theta_1},~\omega_2\in\Omega_{\theta_2}$. Moreover, let $\omega_1\in L_{\theta_1}(0,\infty)$. Then
$$\alpha\ge\bigg(\sum_{i=1}^n\frac{1}{p_{1i}}-\sum_{i=1}^n\frac{1}{q_{1i}}\bigg)_+$$
is necessary for the boundedness of $I_{\alpha}$ from $LM_{\vec{p_1}\theta_1,\omega_1}$ to $LM_{\vec{p_2}\theta_2,\omega_2}$.\\
(2) Let $1\le\vec{p}_1\le\infty,~0<\vec{p}_2\le\infty,~0<\alpha<n,~0<\theta_1,\theta_2<\infty,~\omega_1\in\Omega_{\vec{p}_1\theta_1},~\omega_2\in\Omega_{\vec{p}_2\theta_2}$. Moreover, let $\omega_1\in L_{\theta_1}(0,\infty)$. Then
$$\alpha\ge\bigg(\sum_{i=1}^n\frac{1}{p_{1i}}-\sum_{i=1}^n\frac{1}{q_{1i}}\bigg)_+$$
is necessary for the boundedness of $I_{\alpha}$ from $GM_{\vec{p_1}\theta_1,\omega_1}$ to $GM_{\vec{p_2}\theta_2,\omega_2}$.

\textbf{Proof:} Suppose that
$$\|I_{\alpha}f\|_{LM_{\vec{p}_2\theta_2,\omega_2}}\le\|f\|_{LM_{\vec{p}_1\theta_1,\omega_1}}~~~~\forall f\in LM_{\vec{p}_1\theta_1,\omega_1}.$$
Let $f\in L_{\vec{p}_1}\subseteq LM_{\vec{p}_1\theta_1,\omega_1}$, and $f$ is not almost everywhere 0. Then for $t>1$
$$\|\delta_tf\|_{L_{\vec{p}_{1i}}(Q(0,r))}=t^{-\sum_{i=1}^n\frac{1}{p_{1i}}}\|f\|_{L_{\vec{p}_1}(Q(0,tr))},$$
$$\delta_t^{-1}I_{\alpha}\delta_t=t^{-\alpha}I_{\alpha},$$
$$\|\delta_t^{-1}I_{\alpha}f\|_{L_{\vec{p}_2}(Q(0,r))}=t^{\sum_{i=1}^n\frac{1}{p_{2i}}}\|I_{\alpha}f\|_{L_{\vec{p}_2}(Q(0,r/t))}.$$
So, it is easy for us to acquire that
\begin{eqnarray*}
\|I_{\alpha}f\|_{LM_{\vec{p}_2\theta_2,\omega_2}}&=&t^{\alpha}\|\delta_t^{-1}I_{\alpha}(\delta_t f)\|_{LM_{\vec{p}_2\theta_2,\omega_2}}\\
&=&t^{\alpha+\sum_{i=1}^n\frac{1}{p_{2i}}}\|\omega_2(r)\|I_{\alpha}(\delta_t f)\|_{L_{\vec{p}_2}(Q(0,r/t))}\|_{L_{\theta_2}(0,\infty)}\\
&\le&t^{\alpha+\sum_{i=1}^n\frac{1}{p_{2i}}}\|\omega_2(r)\|I_{\alpha}(\delta_t f)\|_{L_{\vec{p}_2}(Q(0,r))}\|_{L_{\theta_2}(0,\infty)}\\
&\le&Ct^{\alpha+\sum_{i=1}^n\frac{1}{p_{2i}}}\|\omega_1(r)\|\delta_t f\|_{L_{\vec{p}_1}(Q(0,r))}\|_{L_{\theta_1}(0,\infty)}\\
&\le&Ct^{\alpha+\sum_{i=1}^n\frac{1}{p_{2i}}-\sum_{i=1}^n\frac{1}{p_{1i}}}\|\omega_1(r)\|f\|_{L_{\vec{p}_1}(Q(0,tr))}\|_{L_{\theta_1}(0,\infty)}\\
&\le&Ct^{\alpha+\sum_{i=1}^n\frac{1}{p_{2i}}-\sum_{i=1}^n\frac{1}{p_{1i}}}\|\omega_1(r)\|_{L_{\theta_1}(0,\infty)}\|f\|_{L_{\vec{p}_1}(\mathbb{R}^n)}.
\end{eqnarray*}
If $\sum_{i=1}^n\frac{1}{p_{1i}}\le\sum_{i=1}^n\frac{1}{p_{2i}}$, then $0<\alpha<n$ is necessary because of definition of $I_{\alpha}$. If $\sum_{i=1}^n\frac{1}{p_{1i}}>\sum_{i=1}^n\frac{1}{p_{2i}}$, then $\sum_{i=1}^n\frac{1}{p_{1i}}-\sum_{i=1}^n\frac{1}{p_{2i}}<\alpha<n$ is necessary because of the above discussion.

(2) We can prove (2) through the same method as (1), so we omit the proof. $~~~~\blacksquare$

\section{$L_{\vec{p}}\,$-Estimates of $I_{\alpha}$ Over the Cube $Q(x,r)$}\label{sec1}

\par
We consider the following ``partia'' fractional integral operators
$$\underline{I}_{\alpha,r}f(x)=I_{\alpha}(f\chi_{Q(x,r)})(x)=\int_{Q(x,r)}\frac{|f(y)|}{|x-y|^{n-\alpha}}\,dy.$$
$$\bar{I}_{\alpha,r}f(x)=I_{\alpha}(f\chi_{Q(x,r)^{\complement}})(x)=\int_{Q(x,r)^{\complement}}\frac{|f(y)|}{|x-y|^{n-\alpha}}\,dy.$$

\hypertarget{Theorem 4.1}{\textbf{Theorem 4.1.}} Let $0<\vec{p}\le\infty,~0<\alpha<n$ and $f\in L_1^{loc}(\mathbb{R}^n)$. Then for any $Q(x,r)\subseteq\mathbb{R}^n$
$$\|I_{\alpha}(|f|)\|_{L_{\vec{p}}(Q(x,r))}\thicksim\|I_{\alpha}(|f|\chi_{Q(x,2r)})\|_{L_{\vec{p}}(Q(x,r))}+r^{\sum_{i=1}^n\frac{1}{p_i}}\bar{I}_{\alpha,2r}f(x).$$

\textbf{Proof:} Clearly
$$\|I_{\alpha}(|f|)\|_{L_{\vec{p}}(Q(x,r))}\le\|I_{\alpha}(|f|\chi_{Q(x,2r)})\|_{L_{\vec{p}}(Q(x,r))}+\|I_{\alpha}(|f|\chi_{Q(x,2r)^{\complement}})\|_{L_{\vec{p}}(Q(x,r))}.$$
If $y\in Q(x,r),~z\in Q(x,2r)^{\complement}$, then $|y-z|\thicksim|x-z|$. therefor
\begin{eqnarray*}
\|I_{\alpha}(|f|\chi_{Q(x,2r)^{\complement}})\|_{L_{\vec{p}}(Q(x,r))}&=&\|\int_{Q(x,2r)^{\complement}}\frac{|f(z)|}{|y-z|^{n-\alpha}}\,dz\|_{L_{\vec{p}}(Q(x,r))}\\
&\thicksim&\|\int_{Q(x,2r)^{\complement}}\frac{|f(z)|}{|x-z|^{n-\alpha}}\,dz\|_{L_{\vec{p}}(Q(x,r))}\\
&=& \int_{Q(x,2r)^{\complement}}\frac{|f(z)|}{|x-z|^{n-\alpha}}\,dz\|\chi_{Q(x,r)}\|_{L_{\vec{p}}(\mathbb{R}^n)}\\
&\thicksim& r^{\sum_{i=1}^n\frac{1}{p_i}}\bar{I}_{\alpha,2r}f(x).
\end{eqnarray*}

For the left-hand side inequalities, in the one hand,
$$\|I_{\alpha}(|f|\chi_{Q(x,2r)})\|_{L_{\vec{p}}(Q(x,r))}\le\|I_{\alpha}(|f|)\|_{L_{\vec{p}}(Q(x,r))};$$
in the other hand, if $y\in Q(x,r)$ and $z\in Q(x,2r)^{\complement}$, then
\begin{eqnarray*}
\|I_{\alpha}(|f|)\|_{L_{\vec{p}}(Q(x,r))}&\ge&\|I_{\alpha}(|f|\chi_{Q(x,2r)^{\complement}})\|_{L_{\vec{p}}(Q(x,r))}\\
&\thicksim& r^{\sum_{i=1}^n\frac{1}{p_i}}\bar{I}_{\alpha,2r}f(x).
\end{eqnarray*}

The proof is complete. $~~~~\blacksquare$

\hypertarget{Theorem 4.2}{\textbf{Theorem 4.2.}} Let
$$1<\vec{p}_1\le\vec{p}_2<\infty,~\vec{p}_1\neq\vec{p}_2~\sum_{i=1}^n\frac{1}{p_{1i}}-\sum_{i=1}^n\frac{1}{p_{2i}}\le\alpha<n, \hypertarget{4.1}{\eqno{(4.1)}}$$
or
$$0<\vec{p}_2<\vec{p}_1\le\infty,~\vec{p}_1>1,~0<\alpha<n,\hypertarget{4.2}{\eqno{(4.2)}}$$
or
$$0<\vec{p}_2<\infty,~1<\vec{p}_1<\infty,~\sum_{i=1}^n\left(\frac{1}{p_{1i}}-\frac{1}{p_{2i}}\right)_+<\alpha<n, \hypertarget{4.3}{\eqno{(4.3)}}$$
Then
$$\|I_{\alpha}(|f|\chi_{Q(x,2r)})\|_{L_{\vec{p}_2}(Q(x,r))}\lesssim r^{\alpha-\big(\sum_{i=1}^n\frac{1}{p_{1i}}-\sum_{i=1}^n\frac{1}{p_{2i}}\big)}\|f\|_{L_{\vec{p}_1}(Q(x,2r))}.$$

\textbf{Proof:} (1) When $1<\vec{p}_1<\vec{p}_2<\infty$ and
$$\sum_{i=1}^n\frac{1}{p_{1i}}-\sum_{i=1}^n\frac{1}{p_{2i}}\le\alpha<n,$$
it is easy to calculate that
\begin{eqnarray*}
I_{\alpha}(|f|\chi_{Q(x,2r)})(y)&=&\int_{Q(x,2r)}\frac{|f(z)|}{|y-z|^{n-\alpha}}\,dz\\
&\lesssim&r^{\alpha-\big(\sum_{i=1}^n\frac{1}{p_{1i}}-\sum_{i=1}^n\frac{1}{p_{2i}}\big)}\int_{Q(x,2r)}\frac{|f(z)|}{|y-z|^{n-\beta}}\,dz\\
&=&r^{\alpha-\big(\sum_{i=1}^n\frac{1}{p_{1i}}-\sum_{i=1}^n\frac{1}{p_{2i}}\big)}I_{\beta}(f\chi_{Q(x,2r)})(y),
\end{eqnarray*}
where $\beta=\sum_{i=1}^n\frac{1}{p_{1i}}-\sum_{i=1}^n\frac{1}{p_{2i}}$ and $y\in Q(x,r)$. Duce to \hyperlink{Lemma 3.1}{Lemma 3.1}, It is obtained that
\begin{eqnarray*}
\|I_{\alpha}(|f|\chi_{Q(x,2r)})\|_{L_{\vec{p}_2}(Q(x,r))}&\lesssim& r^{\alpha-\big(\sum_{i=1}^n\frac{1}{p_{1i}}-\sum_{i=1}^n\frac{1}{p_{2i}}\big)}\|I_{\beta}(f\chi_{Q(x,2r)})\|_{L_{\vec{p}_2}(Q(x,r))}\\
&\lesssim&r^{\alpha-\big(\sum_{i=1}^n\frac{1}{p_{1i}}-\sum_{i=1}^n\frac{1}{p_{2i}}\big)}\|f\|_{L_{\vec{p}_1}(Q(x,2r))}
\end{eqnarray*}

(2) If $1<\vec{p}_1\le\vec{p}_2<\infty$, we take $\vec{s}=(s_1,s_2,\cdots,s_n)$ such that $1<\vec{s}<\vec{p}_1\le\vec{p}_2<\infty$ and
$$\sum_{i=1}^n\frac{1}{s_i}-\sum_{i=1}^n\frac{1}{p_{2i}}\le\alpha<n.$$
Thanks to (1) and H\"{o}lder's inequality, it is easy to acquire that
\begin{eqnarray*}
\|I_{\alpha}(|f|\chi_{Q(x,2r)})\|_{L_{\vec{p}_2}(Q(x,r))}&\lesssim& r^{\alpha-\big(\sum_{i=1}^n\frac{1}{s_i}-\sum_{i=1}^n\frac{1}{p_{2i}}\big)}\|f\|_{L_{\vec{s}}(Q(x,2r))}\\
&\lesssim&r^{\alpha-\big(\sum_{i=1}^n\frac{1}{p_{1i}}-\sum_{i=1}^n\frac{1}{p_{2i}}\big)}\|f\|_{L_{\vec{p}_1}(Q(x,2r))}.
\end{eqnarray*}

(3) If $0<\vec{p}_2\le\vec{p}_1<\infty,~\vec{p}_1>1$, then
$$0<\alpha<n.$$
By H\"{o}lder's inequality,
\begin{eqnarray*}
\|I_{\alpha}(|f|\chi_{Q(x,2r)})\|_{L_{\vec{p}_2}(Q(x,r))}&\lesssim& r^{\sum_{i=1}^n\frac{1}{p_{2i}}-\sum_{i=1}^n\frac{1}{P_{1i}}}\|I_{\alpha}(|f|\chi_{Q(x,2r)})\|_{L_{\vec{p}_1}(Q(x,r))}.
\end{eqnarray*}
Duce to Minkowski's inequality,
\begin{eqnarray*}
\|I_{\alpha}(|f|\chi_{Q(x,2r)})\|_{L_{\vec{p}_1}(Q(x,r))}&\lesssim& \left\|\int_{Q(x,2r)}\frac{|f(z)|}{|\cdot-z|^{n-\alpha}}\,dz\right\|_{L_{\vec{p}_1}(Q(x,r))}\\
&\le&\left\|\int_{B(0,3\sqrt{n}r)}\frac{|(f\chi_{Q(x,2r)})(\cdot-z)|}{|z|^{n-\alpha}}\,dz\right\|_{L_{\vec{p}_1}(Q(x,r))}\\
&\le& \int_{B(0,3\sqrt{n}r)}\frac{1}{|z|^{n-\alpha}}\,dz \|(f\chi_{Q(x,2r)})\|_{L_{\vec{p}_1}(\mathbb{R}^n)}\\
&\lesssim&r^{\alpha}\|f\|_{L_{\vec{p}_1}(Q(x,2r))}.
\end{eqnarray*}
Hence,
$$\|I_{\alpha}(|f|\chi_{Q(x,2r)})\|_{L_{\vec{p}_2}(Q(x,r))}\lesssim r^{\alpha-(\sum_{i=1}^n\frac{1}{p_{1i}}-\sum_{i=1}^n\frac{1}{p_{2i}})} \|f\|_{L_{\vec{p}_1}(Q(x,2r))}.$$

(4) For other cases that there are $i\neq j$ such that
$$p_{1i}<p_{2i},~p_{1j}>p_{2j}$$
hold, the proof is given as following.

To simplify the process of the proof, let $\vec{p}_1=(\bar{\mathbf{p}}_1,\mathbf{p}'_1),~\bar{\mathbf{p}}_1=(p_{1,1},p_{1,2},\cdots,p_{1,l})$ and $\mathbf{p}'_1=(p_{1,l+1},p_{1,l+2},\cdots,p_{1,n})$. Similarly, let $\vec{p}_2=(\bar{\mathbf{p}}_2,\mathbf{p}'_2),~x=(\bar{x},x'),$ and $z=(\bar{z},z')$. Moreover,
$$\bar{\mathbf{p}}_1\le\bar{\mathbf{p}}_2,~\mathbf{p}'_1>\mathbf{p}'_2.$$
Duce to (1), (2), and (3),
\begin{eqnarray*}
\|I_{\alpha}(|f|\chi_{Q(x,2r)})\|_{L_{\vec{p}_2}(Q(x,r))}&=& \left\|\left\|\int_{Q(x',2r)}\int_{Q(\bar{x},2r)}\frac{|(f\chi_{Q(x,2r)})(\bar{z},z')|}{|\cdot-(\bar{z},z')|^{n-\alpha}}\,d\bar{z}\,dz' \right\|_{L_{\bar{\mathbf{p}}_2}(Q(\bar{x},r))}\right\|_{L_{\mathbf{p}'_2}(Q(x',r))}\\
&\le&\Bigg\|\int_{Q(x',2r)}\frac{1}{|\cdot-z'|^{n-l-\alpha_2}}\\
&\times&\left\|\int_{Q(\bar{x},2r)}\frac{|(f\chi_{Q(x,2r)})(\bar{z},z')|}{|\cdot-\bar{z}|^{l-\alpha_1}} \,d\bar{z}\right\|_{L_{\bar{\mathbf{p}}_2}(Q(\bar{x},r))}\,dz'\Bigg\|_{L_{\mathbf{p}'_2}(Q(x',r))}\\
&\lesssim&r^{\alpha_1-(\sum_{i=1}^l\frac{1}{p_{1i}}-\sum_{i=1}^l\frac{1}{p_{2i}})}\\ &\times&\left\|\int_{Q(x',2r)}\frac{\|(f\chi_{Q(x',2r)})(\cdot,z')\|_{L_{\bar{\mathbf{p}}_1}(Q(\bar{x},2r))}}{|\cdot-z'|^{n-l-\alpha_1}} \,dz'\right\|_{L_{\mathbf{p}'_2}(Q(x',r))}\\
&\lesssim&r^{\alpha-(\sum_{i=1}^n\frac{1}{p_{1i}}-\sum_{i=1}^n\frac{1}{p_{2i}})} \|\|f\|_{L_{\bar{\mathbf{p}}_1}(Q(\bar{x},2r))}\|_{L_{\mathbf{p}'_1}(Q(x',2r))}\\
&\lesssim&r^{\alpha-(\sum_{i=1}^n\frac{1}{p_{1i}}-\sum_{i=1}^n\frac{1}{p_{2i}})} \|f\|_{L_{\vec{p}_1}(Q(x,2r))},
\end{eqnarray*}
where $\alpha=\alpha_1+\alpha_2$ and they satisfy that
$$\sum_{i=1}^l\frac{1}{p_{1i}}-\sum_{i=1}^l\frac{1}{p_{2i}}<\alpha_1<l,~0<\alpha_2<n-l.$$

The proof is complete. $~~~~\blacksquare$

Before starting the next theorem, the following lemma \cite{16} is given.

\hypertarget{Lemma 4.3}{\textbf{Lemma 4.3.}} Let $f$ is be a non-negative measurable function, then for any $r>0$
$$\int_{|x|>r}\frac{f(x)}{|x|^{\beta}}\,dx=\beta\int_r^{\infty}\int_{r\le|x|\le t}f(x)\,dx\frac{dt}{t^{\beta+1}}.$$

\hypertarget{Theorem 4.4}{\textbf{Theorem 4.4.}} Let (\hyperlink{4.1}{4.1}) or (\hyperlink{4.2}{4.2}) or (\hyperlink{4.3}{4.3}) is satisfied. Then
$$\|I_{\alpha}(|f|)\|_{L_{\vec{p}_2}(Q(x,r))}\lesssim r^{\sum_{i=1}^n\frac{1}{p_{2i}}}\int_{r}^{\infty}\|f\|_{L_{\vec{p}_1}(Q(x,t))}\frac{dt}{t^{\sigma+1}},$$
where $\sigma=\sum_{i=1}^n\frac{1}{p_{1i}}-\alpha$.

\textbf{Proof:} Note that if $\sum_{i=1}^n\frac{1}{p_{1i}}\le\alpha$ and $f$ is not almost everywhere equivalent to 0 on $\mathbb{R}^n$, then $\sigma+1\le1$,
$$\int_{r}^{\infty}\|f\|_{L_{\vec{p}_1}(Q(x,t))}\frac{dt}{t^{\sigma+1}}\ge \|f\|_{L_{\vec{p}_1}(Q(x,r))}\int_{r}^{\infty}\frac{dt}{t^{\sigma+1}}=\infty.$$

Let $\alpha<\sum_{i=1}^n\frac{1}{p_{1i}}$. By \hyperlink{Theorem 4.1}{Theorem 4.1} and \hyperlink{Theorem 4.2}{Theorem 4.2}
\begin{eqnarray*}
\|I_{\alpha}(|f|)\|_{L_{\vec{p}_2}(Q(x,r))}&\thicksim&
\|I_{\alpha}(|f|\chi_{Q(x,2r)})\|_{L_{\vec{p}_2}(Q(x,r))}+r^{\sum_{i=1}^n\frac{1}{p_{2i}}}\bar{I}_{\alpha,2r}f(x)\\
&\lesssim&r^{\alpha-\big(\sum_{i=1}^n\frac{1}{p_{1i}}-\sum_{i=1}^n\frac{1}{p_{2i}}\big)}\|f\|_{L_{\vec{p}_1}(Q(x,2r))}+ r^{\sum_{i=1}^n\frac{1}{p_{2i}}}\bar{I}_{\alpha,2r}f(x)\\
&=& \textrm{(I)}+\textrm{(II)}.
\end{eqnarray*}

In the one hand,
\begin{eqnarray*}
\textrm{(I)}&=&r^{\alpha-\big(\sum_{i=1}^n\frac{1}{p_{1i}}-\sum_{i=1}^n\frac{1}{p_{2i}}\big)}\|f\|_{L_{\vec{p}_1}(Q(x,2r))}\\
&\thicksim&r^{\sum_{i=1}^n\frac{1}{p_{2i}}}\int_{2r}^{\infty}\frac{dt}{t^{\sigma+1}}\|f\|_{L_{\vec{p}_1}(Q(x,2r))}\\
&\le&r^{\sum_{i=1}^n\frac{1}{p_{2i}}}\int_{r}^{\infty}\|f\|_{L_{\vec{p}_1}(Q(x,t))}\frac{dt}{t^{\sigma+1}}.
\end{eqnarray*}

In the other hand, thanks to \hyperlink{Lemma 4.3}{Lemma 4.3} and H\"{o}lder's inequality
\begin{eqnarray*}
\textrm{(II)}&=&r^{\sum_{i=1}^n\frac{1}{p_{2i}}}\bar{I}_{\alpha,2r}f(x)\\
&=&r^{\sum_{i=1}^n\frac{1}{p_{2i}}}\int_{Q(x,2r)^{\complement}}\frac{|f(z)|}{|x-z|^{n-\alpha}}\,dz\\
&=&(n-\alpha)r^{\sum_{i=1}^n\frac{1}{p_{2i}}}\int_{2r}^{\infty}\int_{2r\le|x-y|\le t}|f(y)|\,dy\frac{dt}{t^{n-\alpha+1}}\\
&\lesssim&r^{\sum_{i=1}^n\frac{1}{p_{2i}}}\int_{2r}^{\infty}\|f\|_{L_1(Q(x,t))}\frac{dt}{t^{n-\alpha+1}}\\
&\lesssim&r^{\sum_{i=1}^n\frac{1}{p_{2i}}}\int_{r}^{\infty}\|f\|_{L_{\vec{p}_1}(Q(x,t))}\frac{dt}{t^{\sigma+1}}.
\end{eqnarray*}

The proof is complete. $~~~~\blacksquare$

\section{the sufficient condition of the boundedness of $I_{\alpha}$}\label{sec1}

\par
In this section, the boundedness of the fractional Integration is obtained by Hardy operators' boundedness in weighted Lebesgue spaces. According to \cite{17,18}, the boundedness of Hardy operators was obtained. Hence the sufficient condition of the boundedness of $I_{\alpha}$ and its corollaries are given.

Let $H$ be the Hardy operators,
$$(Hg)(t):=\int_0^t g(r)\,dr,~0<t<\infty,$$
where $g(t)$ is non-negative measurable function on $(0,\infty)$.

\hypertarget{Theorem 5.1}{\textbf{Theorem 5.1.}} Let (\hyperlink{4.1}{4.1}) or (\hyperlink{4.2}{4.2}) or (\hyperlink{4.3}{4.3}) is satisfied. Moreover, let $0<\theta_2\le\infty$ and $\omega_2\in \Omega_{\theta_2}$.
Then
$$\|I_{\alpha}f\|_{LM_{\vec{p}_2\theta_2,\omega_2}}\lesssim\|Hg_{\vec{p}_1}\|_{L{\theta_2,\nu_2}(0,,\infty)},$$
where
$$g_{\vec{p}_1}=\|f\|_{L_{\vec{p}_1}(Q(x,t^{-\frac{1}{\sigma}}))},~\sigma=\sum_{i=1}^n\frac{1}{p_{1i}}-\alpha.$$
$$\nu_2(r)=\omega_2(r^{-\frac{1}{\sigma}})r^{-\frac{1}{\sigma}\sum_{i=1}^n\frac{1}{p_{2i}}-\frac{1}{\theta_2\sigma}-\frac{1}{\theta_2}}.$$

\textbf{Proof:} Duce to \hyperlink{Theorem 4.4}{Theorem 4.4},

\begin{eqnarray*}
\|I_{\alpha}f\|_{LM_{\vec{p}_2\theta_2,\omega_2}}&=&\|\omega_2(r)\|I_{\alpha}f\|_{L_{\vec{p}_2}(Q(x,r))}\|_{L_{\theta_2}(0,\infty)}\\
&\lesssim&\|\omega_2(r)r^{\sum_{i=1}^n\frac{1}{p_{2i}}}\int_{r}^{\infty}\|f\|_{L_{\vec{p}_1}(Q(x,t))}\frac{dt}{t^{\sigma+1}}\|_{L_{\theta_2}(0,\infty)}\\
&\thicksim&\|\omega_2(r)r^{\sum_{i=1}^n\frac{1}{p_{2i}}}\int_{0}^{r^{-\frac{1}{\sigma}}}\|f\|_{L_{\vec{p}_1}(Q(x,t^{-\frac{1}{\sigma}}))}dt\|_{L_{\theta_2}(0,\infty)}\\
&\thicksim&\|\omega_2(r^{-\frac{1}{\sigma}})r^{-\frac{1}{\sigma}\sum_{i=1}^n\frac{1}{p_{2i}}-\frac{1}{\theta_2\sigma}-\frac{1}{\theta_2}}Hg_{\vec{p}_1}(r)\|_{L_{\theta_2}(0,\infty)}\\
&=&\|Hg_{\vec{p}_1}\|_{L{\theta_2,\nu_2}(0,,\infty)}.
\end{eqnarray*}

The proof is complete. $~~~~\blacksquare$

\hypertarget{Theorem 5.2}{\textbf{Theorem 5.2.}} Let (\hyperlink{4.1}{4.1}) or (\hyperlink{4.2}{4.2}) or (\hyperlink{4.3}{4.3}) is satisfied. Moreover, let $0<\theta_1,\theta_2\le\infty$, $\omega_1\in \Omega_{\theta_1}$ and $\omega_2\in \Omega_{\theta_2}$. Then if
$$\|Hg_{\vec{p}_1}\|_{L{\theta_2,\nu_2}(0,,\infty)}\lesssim\|g_{\vec{p}_1}\|_{L{\theta_1,\nu_1}(0,,\infty)},$$
where $\upsilon_2$ is same as theorem 5.1 and
$$\nu_1(r)=\omega_1(r^{-\frac{1}{\sigma}})r^{-\frac{1}{\theta_1\sigma}-\frac{1}{\theta_1}},$$
$I_\alpha$ is bounded from $LM_{\vec{p}_1\theta_1,\omega_1}$ to $LM_{\vec{p}_2\theta_2,\omega_2}$.

\textbf{Proof:} Assume that the operator $H$ is boundedness from $L_{\theta_1,\upsilon_1}(0,\infty)$ to $L_{\theta_2,\upsilon_2}(0,\infty)$. By \hyperlink{Theorem 5.1}{Theorem 5.1},
$$\|I_{\alpha}f\|_{LM_{\vec{p}_2\theta_2,\omega_2}}\lesssim\|Hg_{\vec{p}_1}\|_{L{\theta_2,\nu_2}(0,\infty)}\lesssim\|g_{\vec{p}_1}\|_{L_{\theta_1,\nu_1}(0,,\infty)}.$$
Note that
\begin{eqnarray*}
\|g_{\vec{p}_1}\|_{L{\theta_1,\nu_1}(0,,\infty)}&=& \|\nu_1(r)\|f\|_{L_{\vec{p}_1}(Q(x,r^{-\frac{1}{\sigma}})}\|_{L_{\theta_1}(0,\infty)}\\
&=&\|\omega_1(r^{-\frac{1}{\sigma}})r^{-\frac{1}{\theta_1\sigma}-\frac{1}{\theta_1}}\|f\|_{L_{\vec{p}_1}(Q(x,r^{-\frac{1}{\sigma}})}\|_{L_{\theta_1}(0,\infty)}\\
&\thicksim&\|\omega_1(r)\|f\|_{L_{\vec{p}_1}(Q(x,r))}\|_{L_{\theta_1}(0,\infty)}\\
&=&\|f\|_{LM_{\vec{p}_1\theta_1,\omega_1}}.
\end{eqnarray*}

The proof is complete. $~~~~\blacksquare$

According to \hyperlink{Theorem 5.2}{Theorem 5.2} and \cite{17,18}, we get the following Theorem.

\hypertarget{Theorem 5.3}{\textbf{Theorem 5.3.}} Let (\hyperlink{4.1}{4.1}) or (\hyperlink{4.2}{4.2}) or (\hyperlink{4.3}{4.3}) is satisfied. Moreover, let $0<\theta_1,\theta_2\le\infty$, $\omega_1\in \Omega_{\theta_1}$ and $\omega_2\in \Omega_{\theta_2}$.

Then $I_\alpha$ is bounded from $LM_{\vec{p}_1\theta_1,\omega_1}$ to $LM_{\vec{p}_2\theta_2,\omega_2}$, if\\
(1) If $1<\theta_1\le\theta_2<\infty$, then
$$A_1^1:=\sup_{t>0}\bigg(\int_t^{\infty}\omega_2^{\theta_2}(r)r^{\theta_2(\alpha-(\sum_{i=1}^n\frac{1}{p_{1i}}-\sum_{i=1}^n\frac{1}{p_{2i}}))}\,dr\bigg)^{\frac{1}{\theta_2}} \bigg(\int_t^{\infty}\omega_1^{\theta_1}(r)\,dr\bigg)^{-\frac{1}{\theta_1}}<\infty,$$
and
$$A_2^1:=\sup_{t>0}\bigg(\int_0^{t}\omega_2^{\theta_2}(r)r^{\theta_2\sum_{i=1}^n\frac{1}{p_{2i}}}d\,r\bigg)^{\frac{1}{\theta_2}} \bigg(\int_t^{\infty}\frac{\omega_1^{\theta_1}(r)r^{\theta'(\alpha-\sum_{i=1}^n\frac{1}{p_{1i}})}} {\big(\int_r^{\infty}\omega_1^{\theta_1}(\rho)\,d\rho\big)^{-\frac{1}{\theta'_1}}}\,dr\bigg)^{\frac{1}{\theta'_1}}<\infty.$$
(2) If $1<\theta_1\le 1,~0<\theta_1\le\theta_2<\infty$, then $A_1^1<\infty$ and
$$A_2^2:=\sup_{t>0}t^{\alpha-\sum_{i=1}^n\frac{1}{p_{1i}}} \bigg(\int_0^{t}\omega_2^{\theta_2}(r)r^{\theta_2\sum_{i=1}^n\frac{1}{p_{2i}}}d\,r\bigg)^{\frac{1}{\theta_2}} \bigg(\int_t^{\infty}\omega_1^{\theta_1}(r)\,dr\bigg)^{-\frac{1}{\theta_1}}<\infty.$$
(3) If $1<\theta_1<\infty,~0<\theta_2<\theta_1<\infty,~\theta_2\neq 1$, then
\begin{eqnarray*}
A_1^3&:=&\Bigg(\int_0^{\infty}\bigg( \frac{\int_t^{\infty}\omega_2^{\theta_2}(r)r^{\theta_2(\alpha-(\sum_{i=1}^n\frac{1}{p_{1i}}-\sum_{i=1}^n\frac{1}{p_{2i}}))}\,dr} {\int_r^{\infty}\omega_1^{\theta_1}(r)\,dr}\bigg)^{\frac{\theta_2}{\theta_1-\theta_2}} \\
&\times&\omega_2^{\theta_2}(t)t^{\theta_2(\alpha-(\sum_{i=1}^n\frac{1}{p_{1i}}-\sum_{i=1}^n\frac{1}{p_{2i}}))}\Bigg)^{\frac{\theta_1-\theta_2}{\theta_1\theta_2}}<\infty
\end{eqnarray*}
and
\begin{eqnarray*}
A_2^3&:=&\Bigg(\int_0^{\infty}\Bigg[ \big(\int_0^{t}\omega_2^{\theta_2}(r)r^{\theta_2\sum_{i=1}^n\frac{1}{p_{2i}}}d\,r\bigg)^{\frac{1}{\theta_2}} \bigg(\int_t^{\infty}\frac{\omega_1^{\theta_1}(r)r^{\theta'(\alpha-\sum_{i=1}^n\frac{1}{p_{1i}})}} {\big(\int_r^{\infty}\omega_1^{\theta_1}(\rho)\,d\rho\big)^{-\frac{1}{\theta'_1}}}\,dr\bigg)^{\frac{\theta_2-1}{\theta_2}}
\Bigg]^{\frac{\theta_1\theta_2}{\theta_1-\theta_2}}\\
&\times&\frac{\omega_1^{\theta_1}(t)t^{\theta'(\alpha-\sum_{i=1}^n\frac{1}{p_{1i}})}} {\big(\int_t^{\infty}\omega_1^{\theta_1}(\rho)\,d\rho\big)^{-\frac{1}{\theta'_1}}}\,dt\Bigg)^{\frac{\theta_1-\theta_2}{\theta_1\theta_2}} <\infty.
\end{eqnarray*}
(4) If $1=\theta_2<\theta_1<\infty$, then $A_1^3<\infty$ and
\begin{eqnarray*}
A_2^4&:=&\Bigg(\int_0^{\infty}\bigg(\frac{\int_t^{\infty}\omega_2(r)r^{\alpha-(\sum_{i=1}^n\frac{1}{p_{1i}}-\sum_{i=1}^n\frac{1}{p_{2i}})}\,dr +t^{\alpha-\sum_{i=1}^n\frac{1}{p_{1i}}}\int_0^{t}\omega_2^{\theta_2}(r)r^{\theta_2\sum_{i=1}^n\frac{1}{p_{2i}}}\,dr} {\int_t^{\infty}\omega_1^{\theta_1}(r)\,dr}\bigg)^{\theta'_1-1}\\
&\times& t^{\alpha-\sum_{i=1}^n\frac{1}{p_{1i}}}\int_0^{t}\omega_2^{\theta_2}(r)r^{\theta_2\sum_{i=1}^n\frac{1}{p_{2i}}}\,dr
\frac{dt}{t}\Bigg)^{\theta'_1}<\infty.
\end{eqnarray*}
(5) If $0<\theta_2<\theta_1\le 1$, then $A_1^3<\infty$ and
\begin{eqnarray*}
A_2^5&:=&\Bigg(\int_0^{\infty}\sup_{t\le s<\infty}\frac{s^{\alpha-\sum_{i=1}^n\frac{1}{p_{1i}}\frac{\theta_1\theta_2}{\theta_1-\theta_2}}} {\bigg(\int_s^{\infty}\omega_1^{\theta_1}(\rho)\,d\rho\bigg)^{\frac{\theta_2}{\theta_1-\theta_2}}} \bigg(\int_0^{t}\omega_2^{\theta_2}(r)r^{\theta_2\sum_{i=1}^n\frac{1}{p_{2i}}}\,dr\bigg)^{\frac{\theta_2}{\theta_1-\theta_2}}\\
&\times&\omega_2^{\theta_2}(t)t^{\theta_2\sum_{i=1}^n\frac{1}{p_{2i}}}\,dt\Bigg)^{\frac{\theta_1-\theta_2}{\theta_1\theta_2}}<\infty.
\end{eqnarray*}
(6) If $0<\theta_1\le 1,~\theta_2=\infty,$ then
$$A^6:=\mathop{\text{ess sup}}_{0<t\le s<\infty}\frac{\omega_2(2)(t)t^{\sum_{i=1}^n\frac{1}{p_{2i}}}} {s^{\sum_{i=1}^n\frac{1}{p_{1i}}-\alpha}\bigg(\int_s^{\infty}\omega_1^{\theta_1}(r)\,dr\bigg)^{\frac{1}{\theta_1}}}<\infty.$$
(7) If $1<\theta_1<\infty,~\theta_2=\infty$, then
$$A^7:=\mathop{\text{ess sup}}_{t>0}\omega_2(t)t^{\sum_{i=1}^n\frac{1}{p_{2i}}}\Bigg( \int_t^{\infty}\frac{r^{\theta'_1(\alpha-\sum_{i=1}^n\frac{1}{p_{1i}})}} {\bigg(\int_r^{\infty}\omega_1^{\theta_1}(s)\,ds\bigg)^{\theta'_1-1}} \frac{dr}{r}\Bigg)^{\frac{1}{\theta'_1}}<\infty.$$
(8) If $\theta_1=\infty,~0<\theta_2<\infty$, then
\begin{eqnarray*}
A^8&:=&\Bigg(\int_0^{\infty}\bigg(t^{\sum_{i=1}^n\frac{1}{p_{1i}}-\alpha}\int_t^{\infty} \frac{s^{\alpha-\sum_{i=1}^{n}\frac{1}{p_{1i}}-1}\,ds} {\omega_1(s)} \bigg)^{\theta_2}\\
&\times&\omega_2^{\theta_2}(t)t^{\theta_2(\alpha-(\sum_{i=1}^n\frac{1}{p_{1i}}-\sum_{i=1}^n\frac{1}{p_{2i}}))}\,dt \Bigg)^{\frac{1}{\theta_2}}<\infty.
\end{eqnarray*}
(9) If $\theta_1=\theta_2=\infty$, then
$$A^9:=\mathop{\text{ess sup}}_{t>0}\omega_2(t)t^{\sum_{i=1}^n\frac{1}{p_{2i}}} \int_{t}^{\infty}\frac{s^{\alpha-\sum_{i=1}^{n}\frac{1}{p_{1i}}-1}} {\omega_1(s)}\,ds<\infty.$$

Assume that $I_\alpha$ is bounded from $LM_{\vec{p}_1\theta_1,\omega_1}$ to $LM_{\vec{p}_2\theta_2,\omega_2}$. Then
\begin{eqnarray*}
\|I_{\alpha}f\|_{GM_{\vec{p}_2\theta_2,\omega_2}}&=&\sup_{x\in\mathbb{R}^n}\|I_{\alpha}f(\cdot-x)\|_{LM_{\vec{p}_2\theta_2,\omega_2}}\\
&\lesssim&\sup_{x\in\mathbb{R}^n}\|f(\cdot-x)\|_{LM_{\vec{p}_2\theta_2,\omega_2}}\\
&=&\|I_{\alpha}f\|_{GM_{\vec{p}_2\theta_2,\omega_2}}.
\end{eqnarray*}
So we get the following corollary.

\textbf{Corollary 5.4.} If the conditions of \hyperlink{Theorem 5.3}{Theorem 5.3} are satisfied. Moreover, $\omega_1\in\Omega_{\vec{p}_1\theta_1}$ and $\omega_2\in\Omega_{\vec{p}_2\theta_2}$. Then $I_{\alpha}$ are bounded from $GM_{\vec{p}_1\theta_1,\omega_1}$ to $GM_{\vec{p}_2\theta_2,\omega_2}$.

It is well-known that $M_{\alpha}f<I_{\alpha}(|f|)$, then the following corollary is obtained.

\textbf{Corollary 5.5.} If the conditions of \hyperlink{Theorem 5.3}{Theorem 5.3} are satisfied, then $M_\alpha$ are bounded from $LM_{\vec{p}_1\theta_1,\omega_1}$ to $LM_{\vec{p}_2\theta_2,\omega_2}$. Moreover, if $\omega_1\in\Omega_{\vec{p}_1\theta_1}$ and $\omega_2\in\Omega_{\vec{p}_2\theta_2}$. Then $M_{\alpha}$ are also bounded from $GM_{\vec{p}_1\theta_1,\omega_1}$ to $GM_{\vec{p}_2\theta_2,\omega_2}$.

If $\omega_1(r)=r^{\frac{n}{q_1}-\sum_{i=1}^{n}\frac{1}{p_{1i}}},~\omega_2(r)=r^{\frac{n}{q_2}-\sum_{j=1}^{n}\frac{1}{p_{2i}}},$ $\theta_1=\theta_2=\infty$ and , then
$$GM_{\vec{p}_1\theta_1,\omega_1}=\mathcal{M}_{\vec{p}_1}^{q_1}~~GM_{\vec{p}_2\theta_2,\omega_2}=\mathcal{M}_{\vec{p}_2}^{q_2}.$$
Moreover, if $\frac{1}{q_2}=\frac{1}{q_1}-\frac{\alpha}{n}$, then $A^9<\infty$ is satisfied.

\hypertarget{Corollary 5.6}{\textbf{Corollary 5.6.}} Let $\frac{n}{q_1}\le\sum_{i=1}^n\frac{1}{p_{1i}},\frac{n}{q_2}\le\sum_{i=1}^n\frac{1}{p_{2i}}$ and (\hyperlink{4.1}{4.1}) or (\hyperlink{4.2}{4.2}) or (\hyperlink{4.3}{4.3}) is satisfied. Moreover, if
$$\frac{1}{q_2}=\frac{1}{q_1}-\frac{\alpha}{n},$$
then for $f\in\mathcal{M}_{\vec{p}_1}^{q_1}(\mathbb{R}^n)$,
$$\|I_{\alpha}f\|_{\mathcal{M}_{\vec{p}_2}^{q_2}(\mathbb{R}^n)}\le C\|f\|_{\mathcal{M}_{\vec{p}_1}^{q_1}(\mathbb{R}^n)}.$$

\hypertarget{Corollary 5.7}{\textbf{Corollary 5.7.}} Let $0<\alpha<n,~1<\vec{p},\vec{q}<\infty$. Then
$$1<\vec{p}\le\vec{q}<\infty,~\vec{p}\neq\vec{q},~\alpha=\sum_{i=1}^n\frac{1}{p_i}-\sum_{i=1}^n\frac{1}{q_i}.$$
if and only if
$$\|I_{\alpha}f\|_{L_{\vec{q}}}\lesssim\|f\|_{L_{\vec{p}}}$$

\textbf{Proof:} It easy to prove the sufficiency by \hyperlink{Corollary 5.6}{Corollary 5.6}. For the necessity, by \hyperlink{Lemma 3.1}{Lemma 3.1} and \hyperlink{Remark 3.2}{Remark 3.2},
$$\alpha=\sum_{i=1}^n\frac{1}{p_i}-\sum_{i=1}^n\frac{1}{q_i},~1<\vec{p}\le\vec{q}<\infty.$$
It is obvious that $\vec{p}_1\neq\vec{p}_2$. $~~~~\blacksquare$

\end{document}